\documentclass[12pt,eqno,]{article}
\usepackage{amsmath,amsthm,amssymb}

\makeatletter

\newcommand{\Rmnum}[1]{\expandafter\@slowromancap\romannumeral #1@}
\makeatother

\title{Sharp Inequalities between Harmonic, Seiffert, Quadratic and Contraharmonic Means}
\author{\small  Gen-Di Wang$^{1}$, Chen-Yan Yang$^{2}$ and Yu-Ming Chu$^{1}$}
\date{}

\begin{document}
\maketitle

{\footnotesize\rm
\noindent $^{1}$Department of Mathematics, Huzhou Teachers College, Huzhou 313000, China;\\
$^{2}$Department of Mathematics, Zhejiang Sci-Tech University, Hangzhou 310018, China.\\
Correspondence should be addressed to Yu-Ming Chu,
chuyuming@hutc.zj.cn}

\bigskip \noindent{\bf Abstract}: In this paper, we present the
greatest values $\alpha$, $\lambda$ and $p$, and the least values
$\beta$, $\mu$ and $q$ such that the double inequalities $\alpha
D(a,b)+(1-\alpha)H(a,b)<T(a,b)<\beta D(a,b)+(1-\beta) H(a,b)$,
$\lambda D(a,b)+(1-\lambda)H(a,b)<C(a,b)<\mu D(a,b)+(1-\mu) H(a,b)$
and $p D(a,b)+(1-p)H(a,b)<Q(a,b)<q D(a,b)+(1-q)H(a,b)$ hold for all
$a,b>0$ with $a\neq b$, where $H(a,b)=2ab/(a+b)$,
$T(a,b)=(a-b)/[2\arctan((a-b)/(a+b))]$, $Q(a,b)=\sqrt{(a^2+b^2)/2}$,
$C(a,b)=(a^2+b^2)/(a+b)$ and $D(a,b)=(a^3+b^3)/(a^2+b^2)$ are the
harmonic, Seiffert, quadratic, first contraharmonic and second
contraharmonic means of $a$ and $b$, respectively.



\section{Introduction}
\hspace{4mm} For $a,b>0$ with $a\neq b$ the harmonic mean $H(a,b)$,
Seiffert mean $T(a,b)$, quadratic mean $Q(a,b)$, first
contraharmonic mean $C(a,b)$ and second contraharmonic mean $D(a,b)$
are defined by $H(a,b)=2ab/(a+b)$,
$T(a,b)=(a-b)/[2\arctan((a-b)/(a+b))]$, $Q(a,b)=\sqrt{(a^2+b^2)/2}$,
$C(a,b)=(a^2+b^2)/(a+b)$ and $D(a,b)=(a^3+b^3)/(a^2+b^2)$,
respectively. Recently, the harmonic, Seiffert, quadratic, first
contraharmonic and second contraharmonic means have attracted the
attention of many mathematicians. In particular, many remarkable
inequalities for these means can be found in the literature [1-28].
Let $G(a,b)=\sqrt{ab}$, $L(a,b)=(a-b)/(\log{a}-\log{b})$,
$I(a,b)=1/e(b^b/a^a)^{1/(b-a)}$, $A(a,b)=(a+b)/2$, and
$M_{p}(a,b)=[(a^p+b^p)/2]^{1/p}(p\neq 0)$ and $M_{0}(a,b)=\sqrt{ab}$
be the geometric, logarithmic, identric, arithmetic and $p$-th
powers means of $a$ and $b$, respectively. Then it is well known
that the inequalities
\begin{align*}
&\qquad\qquad H(a,b)=M_{-1}(a,b)<G(a,b)=M_{0}(a,b)<L(a,b)<I(a,b)\\
&\quad<A(a,b)=M_{1}(a,b)<T(a,b)<Q(a,b)=M_{2}(a,b)<C(a,b)<D(a,b)
\end{align*}
hold for all $a,b>0$ with $a\neq b$.

For all $a,b>0$ with $a\neq b$ Seiffert [29] established the double
inequality $A(a,b)<T(a,b)<Q(a,b)$. H\"{a}st\"{o} [30] proved that
the function $T(1,x)/M_{p}(1,x)$ is increasing in $[1,\infty)$ if
$p\leq 1$. Chu et al. [31] gave the greatest values $\alpha$ and
$\lambda$, and the least values $\beta$ and $\mu$ such that the
double inequalities $\alpha Q(a,b)+(1-\alpha)A(a,b)<T(a,b)<\beta
Q(a,b)+(1-\beta)A(a,b)$ and
$Q^{\lambda}(a,b)A^{1-\lambda}(a,b)<T(a,b)<Q^{\mu}(a,b)A^{1-\mu}(a,b)$
hold for all $a,b>0$ with $a\neq b$.

For $\alpha_{1},\alpha_{2},\beta_{1},\beta_{2}\in(0,1/2)$, Chu et
al. [32, 33] proved that the double inequalities
\begin{equation*}
Q(\alpha_{1} a+(1-\alpha_{1})b,
\alpha_{1}b+(1-\alpha_{1})a)<T(a,b)<Q(\beta_{1} a+(1-\beta_{1})b,
\beta_{1}b+(1-\beta_{1})a)
\end{equation*}
and
\begin{equation*}
C(\alpha_{2} a+(1-\alpha_{2})b,
\alpha_{2}b+(1-\alpha_{2})a)<T(a,b)<C(\beta_{2} a+(1-\beta_{2})b,
\beta_{2}b+(1-\beta_{2})a)
\end{equation*}
hold for all $a,b>0$ with $a\neq b$ if and only if $\alpha_{1}\leq
(1+\sqrt{16/{\pi}^2-1})/2$, $\beta_{1}\geq (3+\sqrt{6})/6$,
$\alpha_{2}\leq (1+\sqrt{4/\pi-1})/2$ and $\beta_{2}\geq
(3+\sqrt{3})/6$.

In [34] Neuman proved that the double inequalities
\begin{equation*}
\alpha Q(a,b)+(1-\alpha)A(a,b)<NS(a,b)<\beta Q(a,b)+(1-\beta)A(a,b)
\end{equation*}
and
\begin{equation*}
\lambda C(a,b)+(1-\lambda)A(a,b)<NS(a,b)<\mu C(a,b)+(1-\mu)A(a,b)
\end{equation*}
hold for all $a,b>0$ with $a\neq b$ if and only if $\alpha\leq
[1-\log(1+\sqrt{2})]/[(\sqrt{2}-1)\log(1+\sqrt{2})]=0.3249\cdots$,
$\beta \geq 1/3$, $\lambda\leq
[1-\log(1+\sqrt{2})]/\log(1+\sqrt{2})=0.1345\cdots$ and $\beta \geq
1/6$, where
\begin{equation*}
NS(a,b)=\frac{a-b}{2{{\sinh^{-1}}}\left[(a-b)/(a+b)\right]}.
\end{equation*}

It is the aim of this paper to present the greatest values $\alpha$,
$\lambda$ and $p$, and the least values $\beta$, $\mu$ and $q$ such
that the double inequalities
\begin{equation*}
\alpha D(a,b)+(1-\alpha) H(a,b)<T(a,b)<\beta D(a,b)+(1-\beta) H(a,b)
\end{equation*}
\begin{equation*}
\lambda D(a,b)+(1-\lambda)H(a,b)<C(a,b)<\mu D(a,b)+(1-\mu) H(a,b)
\end{equation*}
and
\begin{equation*}
p D(a,b)+(1-p)H(a,b)<Q(a,b)<q D(a,b)+(1-q)H(a,b)
\end{equation*}
hold for all $a,b>0$ with $a\neq b$.

\section{Lemmas}

\hspace{4mm} \setcounter{equation}{0}

In order to prove our main results we need several lemmas, which we
present in this section.

\medskip
{\bf Lemma 2.1.} Let
$f_{1}(t)=({\pi}^2-{\pi}-4)t^6-2\pi(5-\pi)t^5-3\pi(5-\pi)t^4-4(5\pi-{\pi}^2-2)t^3-3\pi(5-\pi)t^2-2\pi(5-\pi)t+({\pi}^2-{\pi}-4)$.
Then there exists $\lambda_{0}>1$ such that $f_{1}(t)<0$ for
$t\in[1,\lambda_{0})$ and $f_{1}(t)>0$ for
$t\in(\lambda_{0},+\infty)$.

\medskip
{\bf Proof.} Simple computations lead to
\begin{equation}
f_{1}(1)=-8\pi(9-2\pi)<0, \quad
\lim\limits_{t\rightarrow+\infty}f_{1}(t)=+\infty,
\end{equation}
\begin{align*}
{f_{1}}'(t)=6(\pi^2-\pi-4)t^5-10\pi(5-\pi)t^4-12\pi(5-\pi)t^3\\
-12(5\pi-{\pi}^2-2)t^2-6\pi(5-\pi)t-2\pi(5-\pi),
\end{align*}
\begin{equation}
{f_{1}}'(1)=-24\pi(9-2\pi)<0, \quad \lim\limits_{t\rightarrow
+\infty}{f_{1}}'(t)=+\infty,
\end{equation}
\begin{equation*}
{f_{1}}''(t)=30({\pi}^2-\pi-4)t^4-40\pi(5-\pi)t^3-36\pi(5-\pi)t^2-24(5\pi-{\pi}^2-2)t-6\pi(5-\pi).
\end{equation*}
\begin{equation}
{f_{1}}''(1)=-8(70\pi+9-17{\pi}^2)<0,\quad \lim\limits_{t\rightarrow
\infty}{f_{1}}''(t)=+\infty,
\end{equation}
\begin{equation*}
{f_{1}}'''(t)=120({\pi}^2-\pi-4)t^3-120\pi(5-\pi)t^2-72\pi(5-\pi)t-24(5\pi-{\pi}^2-2),
\end{equation*}
\begin{equation}
{f_{1}}'''(1)=-48(25\pi+9-7{\pi}^2)<0,\quad
\lim\limits_{t\rightarrow \infty}{f_{1}}'''(t)=+\infty,
\end{equation}
\begin{equation*}
{f_{1}}^{(4)}(t)=360({\pi}^2-\pi-4)t^2-240\pi(5-\pi)t-72\pi(5-\pi),
\end{equation*}
\begin{equation}
{f_{1}}^{(4)}(1)=-96(20\pi+15-7{\pi}^2)<0,\quad
\lim\limits_{t\rightarrow \infty}{f_{1}}^{(4)}(t)=+\infty,
\end{equation}
\begin{align}
{f_{1}}^{(5)}(t)=&720({\pi}^2-\pi-4)t-240\pi(5-\pi)>720({\pi}^2-\pi-4)-240\pi(5-\pi)\nonumber\\
=&960(\pi+1)(\pi-3)>0
\end{align}
for $t>1$.

Inequality (2.6) implies that ${f_{1}}^{(4)}(t)$ is strictly
increasing in $[1,+\infty)$. Then from (2.5) we clearly see that
there exists $\lambda_{1}>1$ such that ${f_{1}}^{(4)}(t)<0$ for
$t\in[1,\lambda_{1})$ and ${f_{1}}^{(4)}(t)>0$ for
$t\in(\lambda_{1},+\infty)$. Hence ${f_{1}}'''(t)$ is strictly
decreasing in $[1,\lambda_{1}]$ and strictly increasing in
$[\lambda_{1},+\infty)$.

It follows from (2.4) and the piecewise monotonicity of
${f_{1}}'''(t)$ that there exists $\lambda_{2}>\lambda_{1}>1$ such
that ${f_{1}}''(t)$ is strictly decreasing in $[1,\lambda_{2}]$ and
strictly increasing in $[\lambda_{2},+\infty)$.

From (2.3) and the piecewise monotonicity of ${f_{1}}''(t)$ we
clearly see that there exists $\lambda_{3}>\lambda_{2}>1$ such that
${f_{1}}'(t)$ is strictly decreasing in $[1,\lambda_{3}]$ and
strictly increasing in $[\lambda_{3},+\infty)$.

The piecewise monotonicity of $f'(t)$ and (2.2) lead to the
conclusion that there exists $\lambda_{4}>\lambda_{3}>1$ such that
${f_{1}}(t)$ is strictly decreasing in $[1,\lambda_{4}]$ and
strictly increasing in $[\lambda_{4},+\infty)$.

Therefore, Lemma 2.1 follows from (2.1) and the piecewise
monotonicity of $f_{1}(t)$.

\medskip
{\bf Lemma 2.2.} Let $c\in(0,1)$, $t>1$ and
\begin{equation}
f_{c}(t)=\frac{t^4-1}{ct^4+(2-c)t^3+(2-c)t+c}-2\arctan\left(\frac{t-1}{t+1}\right).
\end{equation}
Then $f_{4/9}(t)>0$ and $f_{2/\pi}(t)<0$ for all $t>1$.

\medskip
{\bf Proof.} From (2.7) one has
\begin{equation}
f_{c}(1)=0,
\end{equation}
\begin{equation}
\lim\limits_{t\rightarrow +\infty}
f_{c}(t)=\frac{1}{c}-\frac{\pi}{2},
\end{equation}
\begin{equation}
{f_{c}}'(t)=\frac{g_{c}(t)}{(t^2+1)[ct^4+(2-c)t^3+(2-c)t+c]^2},
\end{equation}
where
\begin{align}
g_{c}(t)=&(2-2c^2-c)t^8-4c(2-c)t^7+2c(2-c)t^6+4c^2t^5-2(4c^2-5c+2)t^4\nonumber\\
&+4c^2t^3+2c(2-c)t^2-4c(2-c)t+2-2c^2-c.
\end{align}

We divide the proof into two cases.

{\bf Case 1} $c=4/9$. Then (2.11) becomes
\begin{align}
g_{4/9}(t)=&\frac{2}{81}(47t^8-112t^7+56t^6+32t^5-46t^4+32t^3+56t^2-112t+47)\nonumber\\
=&\frac{2}{81}(t-1)^4(47t^4+76t^3+78t^2+76t+47)>0
\end{align}
for $t>1$.

Therefore, $f_{4/9}(t)>0$ for $t>1$ follows easily from (2.8) and
(2.10) together with (2.12).

{\bf Case 2} $c=2/\pi$. Then (2.9) and (2.11) lead to
\begin{equation}
\lim\limits_{t\rightarrow +\infty}f_{c}(t)=\lim\limits_{t\rightarrow
+\infty}f_{2/\pi}(t)=0,
\end{equation}
\begin{align}
g_{2/\pi}(t)=&\frac{2}{{\pi}^2}[({\pi}^2-\pi-4)t^8-8(\pi-1)t^7+4(\pi-1)t^6+8t^5-2({\pi}^2-5{\pi}+8)t^4\nonumber\\
&+8t^3+4(\pi-1)t^2-8(\pi-1)t+{\pi}^2-\pi-4]\nonumber\\
=&\frac{2(t-1)^2}{{\pi}^2}f_{1}(t),
\end{align}
where $f_{1}(t)$ is defined as in Lemma 2.1.

From (2.1) and (2.14) together with Lemma 2.1 we clearly see that
$f_{2/\pi}(t)$ is strictly decreasing in $[1,\lambda_{0}]$ and
strictly increasing in $[\lambda_{0},+\infty)$.

Therefore, $f_{2/\pi}(t)<0$ for $t>1$ follows from (2.8) and (2.13)
together with the piecewise monotonicity of $f_{2/\pi}(t)$.

\medskip
{\bf Lemma 2.3.} Let $t>1$ and
\begin{equation}
g(t)=\frac{(t^2+1)[(t+1)\sqrt{t^2+1}-2\sqrt{2}t]}{(t^2+t+1)(t-1)^2}.
\end{equation}
Then $g(t)$ is strictly increasing from $(1,\infty)$ onto
$(\sqrt{2}/2,1)$.

\medskip
{\bf Proof.} From (2.15) we get
\begin{equation}
\lim_{t\rightarrow 1^+}g(t)=\frac{\sqrt{2}}{2},\quad
\lim_{t\rightarrow +\infty}g(t)=1.
\end{equation}

Let $t=\tan{x}$, $u=\sin{x}+\cos{x}$. Then $x\in(\pi/4,\pi/2)$,
$u\in(1,\sqrt{2})$ and (2.15) becomes
\begin{align}
g(t)=&\frac{\sin{x}+\cos{x}-\sqrt{2}\sin(2x)}{(1+\sin(2x)/2)(1-\sin(2x))}\nonumber\\
=&\frac{\sin{x}+\cos{x}-\sqrt{2}[(\sin{x}+\cos{x})^2-1]}{\{1+[(\sin{x}+\cos{x})^2-1]/2\}[2-(\sin{x}+\cos{x})^2]}\nonumber\\
=&\frac{2(u-\sqrt{2}u^2+\sqrt{2})}{-u^4+u^2+2}\triangleq h(u).
\end{align}

Equation (2.17) leads to
\begin{align}
h'(u)=&\frac{-4\sqrt{2}u^5+6u^4+8\sqrt{2}u^3-2u^2-12\sqrt{2}u+4}{(-u^4+u^2+2)^2}\nonumber\\
=&\frac{-4\sqrt{2}(u-\sqrt{2})^2}{(-u^4+u^2+2)^2}\left(u^3+\frac{5\sqrt{2}}{4}u^2+u-\frac{\sqrt{2}}{4}\right)<0.
\end{align}

Note that $x\rightarrow\sin{x}+\cos{x}$ is strictly decreasing from
$(\pi/4,\pi/2)$ onto $(1,\sqrt{2})$. Therefore, Lemma 2.3 follows
easily from (2.16)-(2.18).

\section{Main Results}

\hspace{4mm} \setcounter{equation}{0}

{\bf Theorem 3.1.} The double inequality
\begin{equation}
\alpha D(a,b)+(1-\alpha) H(a,b)<T(a,b)<\beta D(a,b)+(1-\beta) H(a,b)
\end{equation}
holds for all $a,b>0$ with $a\neq b$ if and only if $\alpha \leq
4/9$ and $\beta\geq 2/\pi$.

\medskip
{\bf Proof.} Since $H(a,b)$, $T(a,b)$ and $D(a,b)$ are symmetric and
homogeneous of degree 1. Without loss of generality, we assume that
$a>b$. Let $t=a/b>1$ and $c\in(0,1)$. Then simple computations lead
to
\begin{equation}
\frac{T(a,b)-H(a,b)}{D(a,b)-H(a,b)}=\frac{(t^2+1)\left(t^2-4t\arctan\left(\frac{t-1}{t+1}\right)-1\right)}{2(t-1)^2(t^2+t+1)\arctan\left(\frac{t-1}{t+1}\right)},
\end{equation}
\begin{equation}
\lim_{t\rightarrow
1^+}\frac{(t^2+1)\left(t^2-4t\arctan\left(\frac{t-1}{t+1}\right)-1\right)}{2(t-1)^2(t^2+t+1)\arctan\left(\frac{t-1}{t+1}\right)}=\frac{4}{9},
\end{equation}
\begin{equation}
\lim_{t\rightarrow
+\infty}\frac{(t^2+1)\left(t^2-4t\arctan\left(\frac{t-1}{t+1}\right)-1\right)}{2(t-1)^2(t^2+t+1)\arctan\left(\frac{t-1}{t+1}\right)}=\frac{2}{\pi},
\end{equation}
\begin{align}
&T(a,b)-[cD(a,b)+(1-c)H(a,b)]\nonumber\\
&=b\frac{c(t+1)(t^3+1)+2(1-c)t(t^2+1)}{2(t+1)(t^2+1)\arctan\left(\frac{t-1}{t+1}\right)}f_{c}(t),
\end{align}
where $f_{c}(t)$ is defined as in Lemma 2.2.

Therefore, inequality
$4D(a,b)/9+5H(a,b)/9<T(a,b)<2D(a,b)/\pi+(1-2/\pi)H(a,b)$ holds for
all $a,b>0$ with $a\neq b$ follows from (3.5) and Lemma 2.2.

Next, we prove that $\alpha=4/9$ and $\beta=2/\pi$ are the best
possible parameters such that inequality (3.1) holds for all $a,b>0$
with $a\neq b$.

If $\alpha>4/9$, then from (3.2) and (3.3) we know that there exists
$\delta>0$ such that $T(a,b)<\alpha D(a,b)+(1-\alpha) H(a,b)$ for
all $a,b>0$ with $a/b\in(1,1+\delta)$.

If $\beta<2/\pi$, then (3.2) and (3.3) lead to the conclusion that
there exists $T_{0}>1$ such that $T(a,b)>\beta
D(a,b)+(1-\beta)H(a,b)$ for all $a,b>0$ with $a/b\in(T_{0},\infty)$.

\medskip
{\bf Theorem 3.2.} The double inequality
\begin{equation*}
p D(a,b)+(1-p)H(a,b)<Q(a,b)<q D(a,b)+(1-q)H(a,b)
\end{equation*}
holds for all $a,b>0$ with $a\neq b$ if and only if $p\leq 1/2$ and
$q\geq \sqrt{2}/2$.

\medskip
{\bf Proof.} Since $H(a,b)$, $Q(a,b)$ and $D(a,b)$ are symmetric and
homogeneous of degree $1$. Without loss of generality, we assume
that $a>b$. Let $t=a/b>1$, then
\begin{equation}
\frac{Q(a,b)-H(a,b)}{D(a,b)-H(a,b)}=\frac{(t^2+1)[(t+1)\sqrt{t^2+1}-2\sqrt{2}t]}{\sqrt{2}(t^2+t+1)(t-1)^2}=\frac{\sqrt{2}}{2}g(t),
\end{equation}
where $g(t)$ is defined as in Lemma 2.3.

Therefore, Theorem 3.2 follows directly from (3.6) and Lemma 2.3.

\medskip
{\bf Theorem 3.3.} The double inequality
\begin{equation*}
\lambda D(a,b)+(1-\lambda)H(a,b)<C(a,b)<\mu D(a,b)+(1-\mu) H(a,b)
\end{equation*}
holds for all $a,b>0$ with $a\neq b$ if and only if $\lambda\leq
2/3$ and $\mu\geq 1$.

\medskip
{\bf Proof.} Since $H(a,b)$, $C(a,b)$ and $D(a,b)$  are symmetric
and homogeneous of degree $1$. Without loss of generality, we assume
that $a>b$. Let $t=a/b>1$, then simple computations lead to
\begin{equation}
\frac{C(a,b)-H(a,b)}{D(a,b)-H(a,b)}=\frac{t^2+1}{t^2+t+1},
\end{equation}
\begin{equation}
\lim_{t\rightarrow 1^+}\frac{t^2+1}{t^2+t+1}=\frac{2}{3}.
\end{equation}
\begin{equation}
\lim_{t\rightarrow +\infty}\frac{t^2+1}{t^2+t+1}=1.
\end{equation}

Note that the function $t\rightarrow (t^2+1)/(t^2+t+1)$ is strictly
increasing in $[1,\infty)$. Therefore, Theorem 3.3 follows from
(3.7)-(3.9) and the monotonicity of the function $t\rightarrow
(t^2+1)/(t^2+t+1)$.

\bigskip
{\bf Acknowledgements.} This research was supported by the Natural
Science Foundation of China under Grants 11071069 and 11171307, and
the Innovation Team Foundation of the Department of  Education of
Zhejiang Province under Grant T200924.

\end{document}